\documentclass[UTF-8,reqno]{amsart}
\usepackage{enumerate}
\usepackage{amssymb,url,color, booktabs}
\usepackage{mathrsfs}

\usepackage[nobysame]{amsrefs}
\BibSpec{article}{%
+{}{\PrintAuthors} {author}
+{,}{ \textrm} {title}
+{.}{ \textit} {journal}
+{,}{ \textbf} {volume}
+{}{ \parenthesize} {date}
+{,}{ } {pages}
+{.}{ arXiv:} {eprint}
+{.}{} {transition}
}
\BibSpec{book}{%
+{}{\PrintAuthors} {author}
+{,}{ \textit} {title}
+{.}{ \textrm} {series} 
+{,}{ Vol.} {volume} 
+{.}{ } {publisher}
+{,}{ } {date}
+{.}{} {transition}
}
\usepackage{color}
\usepackage[colorlinks=true]{hyperref}
\hypersetup{
    linkcolor=blue,          
    citecolor=red,        
    filecolor=blue,      
    urlcolor=cyan  
}

\numberwithin{equation}{section}

\newcommand{\be}{\begin{eqnarray}}
\newcommand{\ee}{\end{eqnarray}}
\newcommand{\ce}{\begin{eqnarray*}}
\newcommand{\de}{\end{eqnarray*}}
\newtheorem{theorem}{Theorem}[section]
\newtheorem{lemma}[theorem]{Lemma}
\newtheorem{remark}[theorem]{Remark}
\newtheorem{definition}[theorem]{Definition}
\newtheorem{proposition}[theorem]{Proposition}
\newtheorem{Examples}[theorem]{Example}
\newtheorem{corollary}[theorem]{Corollary}

\def\eps{\varepsilon}

\def\e{\mathrm{e}}

\def\[{{\Big[}}
\def\]{{\Big]}}
\def\<{{\langle}}
\def\>{{\rangle}}
\def\({{\Big(}}
\def\){{\Big)}}

\def\bx{{\mathbf{x}}}

\def\dif{{\mathord{{\rm d}}}}

\def\no{\nonumber}
\def\={&\!\!=\!\!&}

\def\cP{{\mathcal P}}

\def\mC{{\mathbb C}}

\def\mE{{\mathbb E}}

\def\mN{{\mathbb N}}

\def\mP{{\mathbb P}}
\def\mQ{{\mathbb Q}}
\def\mR{{\mathbb R}}

\def\1{{\mathbf{1}}}

\def\sA{{\mathscr A}}

\def\sE{{\mathscr E}}
\def\sF{{\mathscr F}}

\def\sI{{\mathscr I}}

\def\sP{{\mathscr P}}

\def\geq{\geqslant}
\def\leq{\leqslant}

\def\eps{\varepsilon}

\def\e{\mathrm{e}}

\def\[{{\Big[}}
\def\]{{\Big]}}
\def\<{{\langle}}
\def\>{{\rangle}}
\def\({{\Big(}}
\def\){{\Big)}}

\def\bx{{\mathbf{x}}}

\def\dif{{\mathord{{\rm d}}}}

\def\no{\nonumber}
\def\={&\!\!=\!\!&}
\def\bt{\begin{theorem}}
\def\et{\end{theorem}}
\def\bl{\begin{lemma}}
\def\el{\end{lemma}}
\def\br{\begin{remark}}
\def\er{\end{remark}}
\def\bx{\begin{Examples}}
\def\ex{\end{Examples}}
\def\bd{\begin{definition}}
\def\ed{\end{definition}}
\def\bp{\begin{proposition}}
\def\ep{\end{proposition}}
\def\bc{\begin{corollary}}
\def\ec{\end{corollary}}

\def\geq{\geqslant}
\def\leq{\leqslant}

\def\<{\langle} \def\>{\rangle}

\def\hh{{h}}

\allowdisplaybreaks

\begin{document}

\title{A discretized version of Krylov's estimate and its applications}


\author{Xicheng Zhang}


\thanks{Research of X. Zhang is partially supported by NNSFC grant of China (No. 11731009)  and the DFG through the CRC 1283  ``Taming uncertainty and profiting from randomness and low regularity in analysis, stochastics and their applications''. }

\address{School of Mathematics and Statistics, Wuhan University, Wuhan, Hubei 430072, P.R.China 
\\ Email: XichengZhang@gmail.com }

\begin{abstract}
In this paper we prove a discretized version of Krylov's estimate for discretized It\^o's processes. As applications, we 
study the weak and strong convergences for Euler's approximation of mean-field SDEs with measurable discontinuous and linear growth coefficients.
Moreover, we also show the propagation of chaos for Euler's approximation of mean-field SDEs.
\end{abstract}

\maketitle

\section{Introduction}

\subsection{Discretized Krylov's estimate}
Let $(\Omega,\sF,\mP; (\sF_t)_{t\geq 0})$ be a complete filtration probability space and $(W_t)_{t\geq 0}$ a $d$-dimensional standard $\sF_t$-Brownian motion. 
Let $\xi_t$ be a $d$-dimensional It\^o process with the following form
\begin{align}\label{Ito}
\xi_t=\xi_0+\int^t_0b_s\dif s+\int^t_0\sigma_s\dif W_s,
\end{align}
where $\xi_0\in\sF_0$, $b_s(\omega):\mR_+\times\Omega\to\mR^d$
and $\sigma_s(\omega): \mR_+\times\Omega\to\mR^d\otimes\mR^d$ are bounded measurable $\sF_t$-adapted processes with bound $\kappa_0$. Suppose that for some $\kappa_1>0$,
$$
\det(\sigma_s(\omega)\sigma_s^*(\omega))\geq\kappa_1,\ \forall (s,\omega)\in\mR_+\times\Omega,
$$
where the asterisk stands for the transpose of a matrix.
It is well known that for any $T>0$ and $p\geq d+1$, there exists a constant $C>0$ depending only on $\kappa_0,\kappa_1,p$ and $d$ 
such that for all $f\in L^p([0,T]\times\mR^{d})$,
\begin{align}\label{Kry9}
\mathbb{E} \left(\int^T_0f(s,\xi_s)\dif s\right)\leq C\| f\|_{L^p([0,T]\times \mR^d)},
\end{align}
and for time-independent $f\in L^p(\mR^d)$ with $p\geq d$,
\begin{align}\label{Kry90}
\mathbb{E} \left(\int^T_0f(\xi_s)\dif s\right)\leq C\| f\|_{L^p(\mR^d)}.
\end{align}
Such estimates were proven by Krylov in \cite{Kr1}, which plays a basic role in the study of SDEs with measurable coefficients (see also \cite{Xi-Zh}
for some extensions).

\medskip

In this paper we are interesting in showing a discretized version of \eqref{Kry9}.  More precisely, for fixed $N\in\mN$, we introduce the following discretized It\^o process:
for $k\in\mN$,
\begin{align}\label{Itos}
\begin{split}
\xi^N_{k}&:=\xi^N_0+\sum_{j=0}^{k-1} b_{j}/N+\sum_{j=0}^{k-1} \sigma_{j}\cdot\big(W_{(j+1)/N}-W_{j/N}\big)\\
&=\xi^N_{k-1}+b_{k-1}/N+\sigma_{k-1}\cdot\big(W_{k/N}-W_{(k-1)/N}\big),
\end{split}
\end{align}
where $\xi^N_0\in\sF_0$, and for each $j\in\mN_0:=\mN\cup\{0\}$, 
$b_j\in \mR^d$ and $\sigma_j\in\mR^d\otimes\mR^d$ are $\sF_{j/N}$-measurable random variables. 
We aim to establish a discretized version of Krylov's estimate for $\xi^N_k$ in the following theorem.
\bt\label{Th0}
Suppose that for some $\kappa_0,\kappa_1>0$ and any $j\in\mN_0$,
$$
|b_j|,\|\sigma_j\|\leq\kappa_0,\ \ \det(\sigma_j\sigma^*_j)\geq\kappa_1,\ \ a.s.
$$
Then for any $p>d+1$, there is a constant $C=C(p,d,\kappa_0,\kappa_1)>0$ such that for any $N\in\mN$ and $f_k\in L^p(\mR^d)$, $k=1,\cdots,N$,
\begin{align}\label{EY1}
\frac{1}{N}\sum_{k=1}^N\mathbb{E}f_k(\xi^N_k)\leq C\left(\frac{1}{N}\sum_{k=1}^{N}\|f_k\|^p_{L^p(\mR^d)}\right)^{1/p}.
\end{align}
Moreover,  for any $p>d$, we have
\begin{align}\label{EY11}
\frac{1}{N}\sum_{k=1}^{N}\mathbb{E}f(\xi^N_k)\leq C\|f\|_{L^p(\mR^d)}.
\end{align}
\et

The motivation of studying the above discretized version of Krylov's estimate comes from the study of Euler's scheme for SDEs with measurable discontinuous coefficients.
Let us consider the following general SDE in $\mR^d$:
\begin{align}\label{SDE00}
\dif X_t=b_t(X_t)\dif t+\sigma_t(X_t)\dif W_t,\ \ X_0=x,
\end{align}
where $b:\mR_+\times\mR^d\to\mR^d$ is a Borel measurable function and $\sigma:\mR_+\times\mR^d\to\mR^d\otimes\mR^d$ is a nondegenerate 
matrix-valued Borel measurable function and continuous in $x$. If $b$ and $\sigma$ are of linear growth in $x$ uniformly in $t$,
it is well known that SDE \eqref{SDE00} admits a unique weak solution $X_t$ (cf. \cite{St-Va}).
Moreover, if in addition $\sigma$ is Lipschitz continuous in $x$ uniformly in $t$,
then SDE \eqref{SDE00} admits a unique strong solution (cf. \cite{Ve}). For $\hh\in(0,1)$,
consider the following Euler approximation of SDE \eqref{SDE00}:
\begin{align}\label{Euler}
\dif X^\hh_t=b_{t_\hh}(X^\hh_{t_\hh})\dif t+\sigma_{t_\hh}(X^\hh_{t_\hh})\dif W_t,\ \ X^\hh_0=x,
\end{align}
where $t_\hh:=[t/\hh]\hh$,
which can be solved recursively as follows:
\begin{align}\label{Euler0}
X^\hh_t=X^\hh_{t_\hh}+b_{t_\hh}(X^\hh_{t_\hh})(t-t_\hh)+\sigma_{t_\hh}(X^\hh_{t_\hh})(W_t-W_{t_\hh}),
\end{align}
or equivalently, for $k=0,1,2,\cdots$ and $t\in[k\hh,(k+1)\hh)$,
\begin{align}\label{Euler1}
X^\hh_t=X^\hh_{k\hh}+b_{kh}(X^\hh_{k\hh})(t-k\hh)+\sigma_{kh}(X^\hh_{k\hh})(W_t-W_{k\hh}).
\end{align}
One would ask whether it holds
\begin{align}\label{DR}
\lim_{\hh\to 0}\mE\left(\sup_{t\in[0,T]}|X^\hh_t-X_t|^2\right)=0,
\end{align}
where the key point of proving the above limit is to show an estimate like \eqref{EY1}. Notice that if we take $\hh=1/N$, then 
$\xi^N_k:=X^{1/N}_{k/N}$ just takes the same form as in \eqref{Itos}. 
In fact,  when $\sigma$ is H\"older continuous, that is, for some $\alpha\in(0,1)$ and $c>0$,
$$
\|\sigma_t(x)-\sigma_t(y)\|\leq c|x-y|^\alpha,
$$
Gy\"ongy and Krylov \cite[Theorem 4.2]{Gy-Kr} proved that $X^\hh_t$ allows a density $\rho^\hh_t(y)$ with 
$$
\left(\int_{\mR^d}|\rho^\hh_t(y)|^q\dif y\right)^{1/q}\leq C(t^{-d/(2p)}+1),\ \ p=\frac{q}{q-1}>\frac{d}{\alpha},
$$
where $C=C(d,p,\kappa_0,\kappa_1)>0$.
From this, it is easy to derive that for any $p>d/\alpha$,
$$
\mE\int^T_0 f(X^\hh_{t_\hh})\dif t\leq C\|f\|_{L^p(\mR^d)}.
$$
The above discretized Krylov estimate plays a key role in \cite{Gy-Kr} for showing \eqref{DR} when $b$ is only  bounded measurable.
However, by \eqref{EY11}, the above estimate holds for any $p>d$ without any continuity assumption on $\sigma$. 
In other words, using \eqref{EY11} we can drop
the {\it continuity} assumption on $\sigma$ in Theorem 2.8 of \cite{Gy-Kr}.
It should be noticed that in the remarkable paper \cite{Gy-Kr}, under very broad assumptions, Gy\"ongy and Krylov used Euler's polygonal approximation to construct 
the strong solution for SDE \eqref{SDE00}.
We mention that if $b$ satisfies some monotonicity condition and $\sigma$ is Lipschitz continuous, Gy\"ongy \cite{Gy} showed 
the rate of almost surely convergence for Euler's scheme. Up to now, there are many works devoted to the study of Euler's approximation for 
SDEs with irregular coefficients under various assumptions, for examples, see \cite{Gi, Le-Sz, Ng-Ta, Ba-Hu-Yu} and references therein.


\subsection{Euler's scheme for DDSDEs}
Another goal of this paper is to use Theorem \ref{Th0} to derive the same results as in \cite{Gy-Kr} for mean-field 
(also called McKean-Vlasov or distribution-dependent in literature) SDEs with measurable discontinuous coefficients 
$b$ and $\sigma$. For $\beta\geq 0$, let $\cP_\beta(\mR^d)$ be the space of all probability measures on $\mR^d$ with finite $\beta$-order moment, 
which is endowed with the weak convergence topology.
Let $\beta\geq 1$. Consider the following distribution-dependent SDE (abbreviated as DDSDE): 
\begin{align}\label{SDE0}
\dif X_t=b_t(X_t,\mu_{X_t})\dif t+\sigma_t(X_t,\mu_{X_t})\dif W_t,\ \ \mbox{Law of $X_0$}=\nu\in\cP_\beta(\mR^d),
\end{align}
where $\mu_{X_t}$ stands for the law of random variable $X_t$ and
$$
b:\mR_+\times\mR^d\times\cP_\beta(\mR^d)\to\mR^d, \ \sigma:\mR_+\times\mR^d\times\cP_\beta(\mR^d)\to\mR^d\otimes\mR^d
$$
are Borel measurable functions. Below we make the following assumptions:
\begin{enumerate}[{\bf (H$_\beta$)}]
\item For each $x,\mu$, $t\mapsto b_t(x,\mu)$ and $\sigma_t(x,\mu)$ are continuous, and for each $t,x$, $\mu\mapsto b_t(x,\mu)$ and $\sigma_t(x,\mu)$ are weakly continuous.
Moreover, for some $\beta\geq 1$, there is a constant $c_0>0$ such that for all $t\geq 0,x\in\mR^d$ and $\mu\in\cP_\beta(\mR^d)$,
$$
|b_t(x,\mu)|+\|\sigma_t(x,\mu)\|\leq c_0(1+|x|+\mu(|\cdot|^\beta)^{1/\beta}),
$$
and the following nondegenerate condition holds: there is a constant $c_1>0$ such that for all $t\geq 0,x\in\mR^d$ and
$\mu\in\cP_\beta(\mR^d)$,
\begin{align}\label{EL}
\det(\sigma\sigma^*)(t,x,\mu)\geq c_1.
\end{align}
\end{enumerate}
\begin{enumerate}[{\bf (H$'_\beta$)}]
\item Let $\bar   b$ and $\bar \sigma$ be two Borel  measurable functions on $\mR_+\times\mR^d\times\mR^d$
with values in $\mR^d$ and $\mR^d\otimes\mR^d$, respectively. Assume that for each $x,y\in\mR^d$,
$t\mapsto \bar  b_t(x,y),\bar\sigma_t(x,y)$ are continuous, and for some $c_0>0$ and all $t\geq 0$, $x,y\in\mR^d$,
$$
|\bar  b_t(x,y)|+\|\bar  \sigma_t(x,y)\|\leq c_0(1+|x|+|y|).
$$
Moreover, for any $\mu\in\cP_\beta(\mR^d)$, define
$$
b_t(x,\mu):=\int_{\mR^d}\bar  b_t(x,y)\mu(\dif y),\ \ \sigma_t(x,\mu):=\int_{\mR^d}\bar  \sigma_t(x,y)\mu(\dif y),
$$
and we also assume the nondegenerate condition \eqref{EL} holds.
\end{enumerate}
The difference between {\bf (H$_\beta$)} and {\bf (H$'_\beta$)} lies that in the later case,
$$
\mu\mapsto b_t(x,\mu), \sigma_t(x,\mu)
$$ may be not continuous with respect to the weak convergence.
Notice that we do not make any continuity assumptions on $\bar  b, \bar\sigma$ in $x,y$.
We now consider the following Euler approximation of DDSDE \eqref{SDE0}:
\begin{align}\label{GR1}
\dif X^\hh_t=b_{t_\hh}\big(X^\hh_{t_\hh},\mu_{X^\hh_{t_\hh}}\big)\dif t+\sigma_{t_\hh}\big(X^\hh_{t_\hh},\mu_{X^\hh_{t_\hh}}\big)\dif W_t,\ \ 
\mbox{Law of $X^h_0$}=\nu.
\end{align}

The following theorem extends \cite[Theorem 2.8]{Gy-Kr} to  DDSDEs.
\bt\label{Th12}
Let $\beta>2$, $\nu\in\cP_\beta(\mR^d)$ and one of {\bf (H$_\beta$)} and {\bf (H$'_\beta$)} holds.
\begin{enumerate}[(i)]
\item Suppose that weak uniqueness holds for DDSDE \eqref{SDE0}. Then there is a unique weak solution $X$ to DDSDE \eqref{SDE0} 
 with initial law $\mP\circ X_0^{-1}=\nu$ so that $X^\hh$ converges to $X$ in distribution.
Moreover, for any bounded measurable $f$,
\begin{align}\label{WE}
\lim_{\hh\to 0}\mE\left(\int^T_0 f(X^\hh_{t_\hh})\dif t\right)=\mE \left(\int^T_0f(X_t)\dif t\right).
\end{align}
\item Suppose that pathwise uniqueness holds for DDSDE \eqref{SDE0}. Then there is a unique strong solution $X$ to DDSDE \eqref{SDE0}  with initial law $\mP\circ X_0^{-1}=\nu$ so that
\begin{align}\label{WE0}
\lim_{\hh\to 0}\mE\left(\sup_{t\in[0,T]}|X^\hh_t-X_t|^2\right)=0.
\end{align}
\end{enumerate}
\et

About the weak and strong uniqueness of DDSDE \eqref{SDE0}, by Girsanov's theorem, 
Li and Min \cite{Li-Mi} obtained the existence and uniqueness of weak solutions 
when $b$ is bounded measurable and $\sigma$ is nondegenerate and Lipschitz continuous.
While under {\bf (H$_\beta$)} or {\bf (H$'_\beta$)}, when $\sigma$
does not depend on $\mu$ and is Lipschitz continuous in $x$ and $b$ is Lipschitz continuous with respect to $\mu$ in case  {\bf (H$_\beta$)}, 
Mishura and Veretennikov \cite{Mi-Ve} showed the strong uniqueness.
In a recent work of the present author with R\"ockner \cite{Ro-Zh}, we established the well-posedness of DDSDEs \eqref{SDE0} with singular drifts
(see also \cite{Hu-Wa}).

\subsection{Propagation of chaos for Euler's scheme}
Below we fix $h\in(0,1)$ and let $\{\xi_j,j\in\mN\}$ be a sequence
of i.i.d. random variables in $\mR^d$ with common distribution $\nu$, and
$\{W^j, j\in\mN\}$ a sequence of independent $d$-dimensional standard Brownian motions. 
For numerical reason, we also consider the following interacting particle approximation for Euler's scheme:
for fixed $N\in\mN$, we define for $j=1,\cdots,N$,
\begin{align}\label{SDE99}
\dif X^{N,j}_t=b_{t_h}\big(X^{N,j}_{t_h},\mu^N_{t_h}\big)\dif t+\sigma_{t_h}\big(X^{N,j}_{t_h},\mu^N_{t_h}\big)\dif W^j_t,\ X^{N,j}_0=\xi_j,
\end{align}
where  $\mu^N_t$ is the empirical measure of $\{X^{N,j}_t, j=1,\cdots,N\}$ defined by
$$
\mu^N_t:=\frac{1}{N}\sum^N_{i=1}\delta_{X^{N,i}_t},
$$
where $\delta_x$ stands for the Dirac measure concentrated at point $x$.
In the following, for simplicity we only consider the case {\bf (H$'_\beta$)}, and in this case  we have
$$
b_{t_h}\big(X^{N,j}_{t_h},\mu^N_{t_h}\big)=\frac{1}{N}\sum_{i=1}^N\bar b_{t_h}\big(X^{N,j}_{t_h},X^{N,i}_{t_h}\big)
$$
and
$$
\sigma_{t_h}\big(X^{N,j}_{t_h},\mu^N_{t_h}\big)=\frac{1}{N}\sum_{i=1}^N\bar\sigma_{t_h}\big(X^{N,j}_{t_h},X^{N,i}_{t_h}\big).
$$
For $j\in\mN$, let $\bar X^j_t$ be the unique solution of the following Euler scheme:
\begin{align}\label{SDE98}
\dif \bar X^{j}_t=b_{t_h}\Big(\bar X^{j}_{t_h}, \mu_{\bar X^j_{t_h}}\Big)\dif t
+\sigma_{t_h}\Big(\bar X^{j}_{t_h}, \mu_{\bar X^j_{t_h}}\Big)\dif W^j_t,\ \ \bar X^j_0=\xi_j.
\end{align}
Clearly, $\{\bar X^j_\cdot, j\in\mN\}$ is  a family of i.i.d. stochastic processes with common distribution as $X^h_\cdot$.
\bt\label{Th13}
Let $\beta>2$ amd $\nu\in\cP_\beta(\mR^d)$.
Suppose that {\bf (H$'_\beta$)} holds and the initial law $\nu$ has a density $\phi\in L^q_{loc}(\mR^d)$ for some $q>1$. 
Then it holds that
for any $T>0$,
\begin{align}\label{DF3}
\lim_{N\to\infty}\sup_{j=1,\cdots,N}\mE\left(\sup_{t\in[0,T]}|X^{N,j}_t-\bar X^j_t|^2\right)=0.
\end{align}
\et
For fixed $h\in(0,1)$ and $N\in\mN$, we use $\sE_h$ and $\sP_N$ to denote the operators of Euler's scheme and 
the interacting particle approximation to DDSDE \eqref{SDE0}, respectively:
$$
\sE_h: X\to X^h,\ \ \sP_N: X\to (X^{N,j})_{j=1,\cdots,N}.
$$
From the construction, it is easy to see that
$$
\sE_h\sP_N=\sP_N\sE_h.
$$
Under {\bf (H$'_\beta$)}, suppose that the pathwise uniqueness holds for DDSDE \eqref{SDE0}. 
Then by Theorems \ref{Th12} and \ref{Th13}, we have
$$
\lim_{h\to 0}\lim_{N\to\infty}\mE\|\sP_N \sE_h X-X\|^2_{C([0,T])}=0.
$$
Here an open question is to show that
\begin{align}\label{GFL}
\lim_{N\to\infty}\lim_{h\to 0}\mE\|\sP_N \sE_h X-X\|^2_{C([0,T])}=0
\end{align}
and
$$
\lim_{N\to\infty}\mE\|\sP_N \sE_{1/N} X-X\|^2_{C([0,T])}=0.
$$
Obviously, the obstacle is to show the following propagation of chaos under {\bf (H$'_\beta$)}:
\begin{align}\label{GFH}
\lim_{N\to\infty}\mE\|\sP_N X-X\|^2_{C([0,T])}=0.
\end{align}
When $b$ and $\sigma$ are Lipschitz continuous in $x$ and $\mu$, the above propagation of chaos \eqref{GFH} was proven by Sznitman \cite{Sz}.
Recently, Bao and Huang \cite{Ba-Hu} proved \eqref{GFL} by Zvonkin's transformation 
when $b$ and $\sigma$ are H\"older continuous in $x$ and Lipschitz continuous in $\mu$
with respect to the Wasserstein distance. However, under {\bf (H$'_\beta$)}, proving \eqref{GFH} seems to be a challenge problem.

\subsection{Plan and Notations}
This paper is organized as follows: In Section 2, we prove Theorem \ref{Th0}. In Section 3, we prove Theorem \ref{Th12}.
In Section 4, we prove Theorem \ref{Th13}.
Throughout this paper we use the following conventions:
\begin{itemize}
\item For a matrix $\sigma$, we use $\|\sigma\|$ to denote the Hilbert-Schmidt norm of $\sigma$.
\item For $R>0$, we use $B_R$ to denote the ball in $\mR^d$ with radius $R$ and center $0$.
\item We use $A\lesssim B$ (resp. $\asymp$) to denote $A\leq CB$ (resp. $C^{-1}B\leq A\leq CB$) for some unimportant constant $C\geq 1$, 
whose dependence on the parameters can be traced from the context. 
\end{itemize}

\section{Proof of Theorem \ref{Th0}}

To prove \eqref{EY1}, we shall use the classical Krylov estimate \eqref{Kry9}.
For this we need to embed $\xi^N_k$ into a continuous It\^o process. For $k\in\mN_0$ and $t\in [k/N, (k+1)/N)$, we define
$$
\tilde b^N_t:=b_k,\ \ \tilde \sigma^N_t:=\sigma_k
$$
and
\begin{align*}
X^N_t&:=\xi^N_k+\tilde b^N_t\cdot (t-k/N)+\tilde \sigma^N_t\cdot \big(W_t-W_{k/N}\big).
\end{align*}
In this way, it is easy to see that $X^N_{k/N}=\xi^N_k$ and
$$
X^N_t=\xi^N_0+\int^t_0\tilde b^N_s\dif s+\int^t_0\tilde \sigma^N_s\dif W_s,\ \ t\geq 0.
$$
Similarly, let $(f_k)_{k\in\mN}$ be a family of nonnegative measurable functions in $\mR^d$. If we define
$$
\tilde f_N(t,x):=\sum_{k=0}^\infty f_{k+1}(x)\1_{t\in[k/N, (k+1)/N)},\ \ t\geq 0, x\in\mR^d,
$$
then we can write
$$
\frac{1}{N}\sum_{k=1}^N\mathbb{E}f_k(\xi^N_k)
=\int^1_0\mathbb{E}\tilde f_N\big(t,X^N_{t^N_+}\big)\dif t,
$$
where $t^N_+:=([tN]+1)/N$. Moreover, by \eqref{Kry9} we have for any $p\geq d+1$ and $y\in\mR^d$,
\begin{align}\label{Kry10}
\int^1_0\mathbb{E}\left(\tilde f_N(t,X^N_t+y)\right)\dif t\leq C\|\tilde f_N\|_{L^p([0,1]\times\mR^d)}
=C\left(\frac{1}{N}\sum_{k=1}^{N}\|f_k\|^p_p\right)^{1/p}.
\end{align}
At this moment, we can not immediately conclude \eqref{EY1} because we need to treat $\mE\tilde f_N\big(t,X^N_{t^N_+}\big)$ rather than $\mE\tilde f_N(t,X^N_t)$.
Notice that for $t\in[(k-1)/N,k/N)$,
\begin{align*}
\xi^N_k=X^N_{t}+b_{k-1}\cdot \ell^k_t+\sigma_{k-1}\cdot \big(W_{k/N}-W_t\big),
\end{align*}
where
$$
\ell^k_t:=k/N-t.
$$
Since $W_{k/N}-W_t$ is independent with $X^N_t$ and $b_{k-1}, \sigma_{k-1}$, by the change of variable we have
\begin{align*}
\frac{1}{N}\mathbb{E}f_{k}(\xi^N_{k})
&=\mE\int^{k/N}_{(k-1)/N}\! f_k\big(X^N_t+b_{k-1}\ell^k_t+\sigma_{k-1}\cdot(W_{k/N}-W_t)\big)\dif t\\
&=\int^{k/N}_{(k-1)/N}\!\int_{\mR^d}\mE f_k\big(X^N_t+b_{k-1}\ell^k_t+\sigma_{k-1}y\big)\varphi_{\ell^k_t}(y)\dif y\dif t\\
&=\int^{k/N}_{(k-1)/N}\!\mE\int_{\mR^d}f_k\big(X^N_t+y\big)
\varphi_{\ell^k_t}\big(\sigma^{-1}_{k-1}(y-b_{k-1}\ell^k_t)\big)|\det\sigma_{k-1}^{-1}|\dif y\dif t,
\end{align*}
where $\varphi_t(y):=(2\pi t)^{-d/2}\e^{-\frac{|y|^2}{2t}}$ is the distributional density of Brownian motion $W_t$. Noticing that
$$
|\det \sigma_{k-1}^{-1}|=1/|\det \sigma_{k-1}|\leq 1/\sqrt{\kappa_1},
$$
and
\begin{align*}
|\sigma^{-1}_{k-1}(y-b_{k-1}\ell^k_t)|^2&\geq c_0|y-b_{k-1}\ell^k_t|^2\geq \tfrac{c_0}{2}\big(|y|^2-2|b_{k-1}|^2|\ell^k_t|^2\big),
\end{align*}
we have for $\lambda=\sqrt{c_0/2}$ and $t\in[(k-1)/N,k/N)$,
\begin{align*}
\varphi_{\ell^k_t}(\sigma^{-1}_{k-1}(y-b_{k-1}\ell^k_s))\leq\e^{c_0\kappa_0^2|\ell^k_t|/2} \varphi_{\ell^k_t}(\lambda y)
\leq\e^{c_0\kappa_0^2/(2N)} \varphi_{\ell^k_t}(\lambda y).
\end{align*}
Hence, for $\gamma=\frac{p}{d+1}$ and $q=\frac{\gamma}{\gamma-1}$, by H\"older's inequality we further have
\begin{align*}
\frac{1}{N}\sum_{k=1}^{N}\mathbb{E}f_{k}(\xi^N_{k})&\lesssim \int_{\mR^d}\left(\sum_{k=1}^{N}\int^{k/N}_{(k-1)/N}\mE f_k\big(X^N_t+y\big)
\varphi_{\ell^k_t}(\lambda y)\dif t\right)\dif y\\
&\leq\int_{\mR^d}\!\!\left(\sum_{k=1}^{N}\int^{k/N}_{(k-1)/N}\mE f_k\big(X^N_t+y\big)^\gamma\dif t\right)^{\frac{1}{\gamma}}\!\!\!
\left(\sum_{k=1}^{N}\int^{k/N}_{(k-1)/N}|\varphi_{\ell^k_t}(\lambda y)|^q\dif t\right)^{\frac{1}{q}}\!\!\dif y\\
&=\int_{\mR^d}\left(\mE\int^1_0\tilde f_N(t,X^N_t+y)^\gamma\dif t\right)^{\frac{1}{\gamma}}
\left(N\int^{1/N}_0|\varphi_t(\lambda y)|^q\dif t\right)^{\frac{1}{q}}\dif y.
\end{align*}
By Krylov's estimate \eqref{Kry10}, we obtain
\begin{align}\label{EY10}
\frac{1}{N}\sum_{k=1}^{N}\mathbb{E}f_{k}(\xi^N_{k})
\lesssim \|\tilde f_N\|_{L^{\gamma(d+1)}([0,1]\times\mR^d)}\int_{\mR^d}\left(N\int^{1/N}_0|\varphi_t(\lambda y)|^q\dif t\right)^{\frac{1}{q}}\dif y.
\end{align}
Notice that by the change of variable and the scaling property of $\varphi_t(y)$,
\begin{align*}
&\int_{\mR^d}\left(N\int^{1/N}_0\varphi_t(\lambda y)^q\dif t\right)^{\frac{1}{q}}\dif y
=\int_{\mR^d}\left(\int^1_0\varphi_{t/N}(\lambda y)^q\dif t\right)^{\frac{1}{q}}\dif y\\
&\quad=N^{d/2}\int_{\mR^d}\left(\int^1_0\varphi_t(\lambda y\sqrt{N})^q\dif t\right)^{\frac{1}{q}}\dif y
=\lambda^{-d}\int_{\mR^d}\left(\int^1_0\varphi_t(y)^q\dif t\right)^{\frac{1}{q}}\dif y.
\end{align*}
The desired estimate \eqref{EY1} now follows by \eqref{EY10} and showing that the last integral is finite. In fact,
by the change of variable, for some $c=c(d,q)>0$,
\begin{align*}
&\int_{\mR^d}\left(\int^1_0\varphi_t(y)^q\dif t\right)^{\frac{1}{q}}\dif y
=c\int_{\mR^d}|y|^{\frac{2}{q}-d}\left(\int^\infty_{|y|^2}t^{\frac{dq}{2}-2}\e^{-t}\dif t\right)^{\frac{1}{q}}\dif y\\
&\lesssim\int_{B_1}|y|^{\frac{2}{q}-d}\left(\int^\infty_0t^{\frac{dq}{2}-2}\e^{-t}\dif t\right)^{\frac{1}{q}}\dif y
+\int_{B^c_1}|y|^{\frac{2}{q}-d}\left(\int^\infty_{|y|^2}t^{\frac{dq}{2}-2}\e^{-t}\dif t\right)^{\frac{1}{q}}\dif y<\infty.
\end{align*}
As for \eqref{EY11} it follows by using \eqref{Kry90} in the above proof.

\br
When $b_j$ and $\sigma_j$ are nonrandom, the estimate \eqref{EY1} is trivial
because $\xi^N_k$, $k=1,\cdots,N$ are Gaussian random variables. However, 
in the general case, we only know that $\xi^N_k$ is a nondegenerate semimartingale with respect to $\sF_{k/N}$.
Here it is quite interesting to give a purely probabilistic proof for Theorem \ref{Th0}.
It should be noticed that \eqref{Kry9} can be derived from \eqref{EY1} by discretized approximation.
\er

\section{Proof of Theorem \ref{Th12}}

The following lemma is standard by Burkholder and Gronwall's inequalities.
\bl\label{Le31}
Let $\beta>2$. Under {\bf (H$_\beta$)} or {\bf (H$'_\beta$)}, for any $T>0$, there is a constant $C>0$ such that for all $h\in(0,1)$,
\begin{align}\label{KJ1}
\mE\left(\sup_{t\in[0,T]}|X_t^h|^\beta\right)\leq C(1+\mE|X_0|^\beta),
\end{align}
and for any $s,t\in[0,T]$,
\begin{align}\label{KJ4}
\mE|X^h_s-X^h_t|^\beta\leq C|s-t|^{\beta/2}.
\end{align}
\el
\begin{proof}
Notice  that
\begin{align}\label{KL1}
X^h_t=X_0+\int^t_0 b_{s_h}\big( X^h_{s_h},\mu_{X^h_{s_h}}\big)\dif s
+\int^t_0\sigma_{s_h}\big( X^h_{s_h},\mu_{X^h_{s_h}}\big)\dif W_s.
\end{align}
For simplicity, we let $|X^h_t|_*:=\sup_{s\in[0,t]}|X^h_s|$.
By Burkholder's inequality and  the linear growth of $b$ and $\sigma$, we have
\begin{align*}
\mE|X^h_t|_*^\beta&\lesssim \mE|X_0|^\beta+\int^t_0 \Big(1+\mE|X^h_{s_h}|^\beta+\mu_{X^h_{s_h}}(|\cdot|^\beta)\Big)\dif s\\
&\lesssim \mE|X_0|^\beta+\int^t_0 \Big(1+\mE |X^h_s|_*^\beta\Big)\dif s,
\end{align*}
which implies \eqref{KJ1} by Gronwall's inequality.
As for \eqref{KJ4}, it follows by \eqref{KL1} and \eqref{KJ1}.
\end{proof}

Let $\mQ_\hh$ be the law of $(X^\hh_\cdot, W_\cdot)$ in product space $\mC\times\mC$, where $\mC$ is the continuous functions space. 
By \eqref{KJ4}, since $\beta>2$,
$(\mQ_\hh)_{\hh\in(0,1)}$ is tight. Therefore, by Prokhorov's theorem, there are a subsequence $h_n\to 0$ as $n\to\infty$ and $\mQ\in\cP(\mC\times\mC)$ so that
$$
\mQ_n:=\mQ_{h_n}\to \mQ\mbox{ weakly.}
$$
Now, by Skorokhod's representation theorem,  there are a probability space $(\tilde\Omega,\tilde\sF,\tilde\mP)$ and 
random variables $(\tilde X^n , \tilde W^n)$ and $(\tilde X,\tilde W)$ defined on it such that
\begin{align}\label{FD4}
(\tilde X^n , \tilde W^n)\to (\tilde X,\tilde W),\ \ \tilde\mP-a.s.
\end{align}
and
\begin{align}\label{FD5}
\tilde\mP\circ(\tilde X^n , \tilde W^n)^{-1}=\mQ_n=\mP\circ(X^{h_n} , W)^{-1},\quad
\tilde\mP\circ(\tilde X, \tilde W)^{-1}=\mQ.
\end{align}
Define $\tilde\sF^n _t:=\sigma(\tilde W^n _s, \tilde X^n_s;s\leq t)$. 
Notice that
\begin{align*}
&\mP(W _t-W_s\in\cdot |\sF _s)=\mP(W _t-W _s\in\cdot)\Rightarrow
\tilde\mP(\tilde W^n _t-\tilde W^n _s\in\cdot |\tilde \sF^n _s)=\tilde\mP(\tilde W^n _t-\tilde W^n _s\in\cdot).
\end{align*}
In other words, $\tilde W^n_t$ is an $\tilde\sF_t^n $-Brownian motion. Thus, by \eqref{KL1} and \eqref{FD5} we have
\begin{align}\label{KK1}
\tilde X^n _t=\tilde X^n_0+\int^t_0b_{s_n}\big(\tilde X^n_{s_n},\mu_{\tilde X^n_{s_n}}\big)\dif s
+\int^t_0\sigma_{s_n}\big(\tilde X^n_{s_n},\mu_{\tilde X^n_{s_n}}\big)\dif \tilde W^n_s,
\end{align}
where 
$s_n:=s_{h_n}=[s/h_n]h_n.$

To take the limits, we recall a result of Skorokhod \cite[p.32]{Sk}. 
\bl\label{Le332}
Let $\{f_n(t),t\geq 0, n\in\mN\}$ be a sequence of measurable $\tilde\sF^n_t$-adapted processes. Suppose that
\begin{enumerate}[(i)]
\item For every $T, \eps>0$, there is an $M_\eps>0$ such that for all $n$,
$$
\tilde\mP\left\{\sup_{t\in[0,T]}|f_n(t)|>M_\eps\right\}\leq\eps.
$$
\item For each $t$, $f_n(t)\to f(t)$ in probability as $n\to\infty$, and for every $T, \eps>0$, 
$$
\lim_{\delta\to 0}\lim_{n\to\infty}\sup_{|t-s|\leq\delta, s,t\in[0,T]}\tilde\mP(|f_n(t)-f_n(s)|>\eps)=0,
$$
or for every $T,\eps>0$,
$$
\lim_{n\to\infty}\tilde\mP\left\{\sup_{t\in[0,T]}|f_n(t)-f(t)|>\eps\right\}=0.
$$
\end{enumerate}
Then it holds that for every $T>0$, 
$$
\int^T_0 f_n(t)\dif \tilde W^n_t\stackrel{n\to\infty}{\to}\int^T_0 f(t)\dif \tilde W_t,\mbox{ in probability}.
$$
\el

Using the above lemma we can show the following limits by the discretized Krylov estimate.

\bl\label{Le32}
Under {\bf (H$_\beta$)} or {\bf (H$'_\beta$)},  for each $t>0$, the following limits hold
\begin{align}
\int^t_0b_{s_n}\left( \tilde X^n _{s_n},\mu_{\tilde X^n_{s_n}}\right)\dif s
&\to \int^t_0b\left(s,\tilde X_{s},\mu_{\tilde X_{s}}\right)\dif s,\label{Lim1}\\
\int^t_0\sigma_{s_n}\left( \tilde X^n _{s_n},\mu_{\tilde X^n_{s_n}}\right)\dif \tilde W^n _s
&\to \int^t_0\sigma\left(s,\tilde X_{s},\mu_{\tilde X_{s}}\right)\dif \tilde W_s\label{Lim2}
\end{align}
in probability as $n\to 0$,
\el
\begin{proof}
We only prove \eqref{Lim2} in case {\bf (H$'_\beta$)}. The others are similar and easier. 
Below for simplicity we shall drop the tilde.  In case {\bf (H$'_\beta$)}, define
$$
\bar\sigma^\eps_t(x,y):=\bar\sigma_t(\cdot,\cdot)*\varrho_\eps(x,y),\ \ \sigma^\eps_t(x,\mu):=\int_{\mR^d}\bar\sigma^\eps_t(x,y)\mu(\dif y),
$$
where $(\varrho_\eps)_{\eps\in(0,1)}$ is a family of mollifiers in $\mR^d\times\mR^d$ with support in $B_\eps\times B_\eps$.
For fixed $\eps\in(0,1)$, since $\bar\sigma_\eps$ is continuous and linear growth in $x,y$, by \eqref{FD4} and Lemma \ref{Le332} (see also \cite[Lemma 3.1]{Gy-Kr}),
it is easy to see that for fixed $\eps\in(0,1)$,
$$
\int^t_0\sigma^\eps_{s_n}\left( X^n _{s_n},\mu_{ X^n_{s_n}}\right)\dif  W^n _s
\to \int^t_0\sigma^\eps_s\left(X_{s},\mu_{ X_{s}}\right)\dif  W_s
$$
in probability as $n\to\infty$. Indeed, it suffices to prove the following two limits:
\begin{align}
\int^t_0\sigma^\eps_{s_n}\left( X^n _{s_n},\mu_{ X^n_{s_n}}\right)\dif  W^n _s&\to \int^t_0\sigma^\eps_s\left(X^n_{s},\mu_{ X^n_{s}}\right)\dif  W^n_s,\label{Es1}\\
\int^t_0\sigma^\eps_s\left( X^n _{s},\mu_{ X^n_{s}}\right)\dif  W^n _s&\to \int^t_0\sigma^\eps_s\left(X_{s},\mu_{ X_{s}}\right)\dif  W_s\label{Es2}
\end{align}
in probability as $n\to\infty$. Limit \eqref{Es1} follows by \eqref{KJ4} and the continuity of $t\mapsto\sigma^\eps_t(x,y)$,
and limit \eqref{Es2} follows by \eqref{FD4} and Lemma \ref{Le332}.
Therefore, it remains to prove that
\begin{align}\label{KK3}
\int^t_0\sigma^\eps_{s_n}\left( X^n _{s_n},\mu_{ X^n_{s_n}}\right)\dif  W^n _s
\to \int^t_0\sigma_{s_n}\left( X^n _{s_n},\mu_{ X^n_{s_n}}\right)\dif  W^n _s
\end{align}
in probability uniformly in $n$ as $\eps\to0$, and
\begin{align}\label{KK4}
\int^t_0\sigma^\eps_s\left( X_{s},\mu_{X_{s}}\right)\dif  W_s
\to \int^t_0\sigma_{s}\left( X_{s},\mu_{ X_{s}}\right)\dif  W_s\mbox{ in probability as $\eps\to0$}.
\end{align}
We only show \eqref{KK3}. By It\^o's isometric formula, we have
\begin{align*}
&\mE\left|\int^t_0\left[\sigma^\eps_{s_n}\left( X^n _{s_n},\mu_{ X^n_{s_n}}\right)-\sigma_{s_n}\left( X^n _{s_n},\mu_{ X^n_{s_n}}\right)\right]\dif  W^n _s\right|^2\\
&=\int^t_0\mE\left\|\sigma^\eps_{s_n}\left( X^n _{s_n},\mu_{ X^n_{s_n}}\right)-\sigma_{s_n}\left( X^n _{s_n},\mu_{ X^n_{s_n}}\right)\right\|^2\dif s\\
&\leq\int^t_0\mE\left\|\bar\sigma^\eps_{s_n}\left( X^n _{s_n},\bar X^n_{s_n}\right)-\bar\sigma_{s_n}\left( X^n _{s_n},\bar X^n_{s_n}\right)\right\|^2\dif s=:J^n_\eps(t),
\end{align*}
where $\bar X^n_\cdot$ is an independent copy of $X^n_\cdot$. More precisely,  $(X^n,\bar X^n)$ solves the following equation (Euler scheme):
\begin{align}\label{EU21}
\left\{
\begin{aligned}
\dif  X^{n}_t=b_{t_n}\big(X^{n}_{t_n},\mu_{ X^n_{t_n}}\big)\dif t+\sigma_{t_n}\big(X^{n}_{t_n},\mu_{ X^n_{t_n}}\big)\dif W_t,\\
\dif  \bar X^{n}_t=b_{t_n}\big(\bar X^{n}_{t_n},\mu_{ \bar X^n_{t_n}}\big)\dif t+\sigma_{t_n}\big(\bar X^{n}_{t_n},\mu_{ \bar X^n_{t_n}}\big)\dif\bar W_t,
\end{aligned}
\right.
\end{align}
where $(W,X^n_0)$ and $(\bar W, \bar X^n_0)$ are independent and have the same distributions.
In order to use the discretized Krylov estimate to show
\begin{align}\label{KQ1}
\lim_{\eps\to 0}\sup_n J^n_\eps(t)=0,
\end{align}
we use a standard stopping time technique. For $R>0$, we define a stopping time
$$
\tau^n_R:=\inf\{t>0: |X^n_t|\vee|\bar X^n_t|>R\},
$$
and make the following decomposition:
\begin{align*}
J^n_\eps(t)&=\int^t_0\mE\Big(\1_{t\geq\tau^n_R}\|\bar\sigma_{s_n}(  X^n _{s_n},\bar X^n_{s_n})
-\bar\sigma^\eps_{s_n}( X^n_{s_n},\bar X^n_{s_n})\|^2\Big)\dif s\\
&+\int^t_0\mE\Big(\1_{t<\tau^n_R}\|\bar\sigma_{s_n}(  X^n _{s_n},\bar X^n_{s_n})-\bar\sigma^\eps_{s_n}( X^n_{s_n},\bar X^n_{s_n})\|^2\Big)\dif s=:J^{n,1}_{R,\eps}(t)+J^{n,2}_{R,\eps}(t).
\end{align*}
For $J^{n,1}_{R,\eps}(t)$, by H\"older's inequality, \eqref{KJ1} and Chebyshev's inequality we have
\begin{align}
J^{n,1}_{R,\eps}(t)&\leq\mP(t\geq\tau^n_R)^{\frac{\beta-2}{\beta}}\left(\int^t_0\mE\|\bar\sigma_{s_n}(  X^n _{s_n},\bar X^n_{s_n})
-\bar\sigma^\eps_{s_n}( X^n_{s_n},\bar X^n_{s_n})\|^\beta\dif s\right)^{\frac{2}{\beta}}\no\\
&\lesssim\left(\frac{\mE\big(\sup_{t\in[0,T]}|X^n_t|\vee|\bar X^n_t|^\beta\big)}{R^\beta}\right)^{\frac{\beta-2}{\beta}}\!\!\!\left(\int^t_0
\left(1+\mE|X^n _{s_n}|^\beta+\mE|\bar X^n_{s_n}|^\beta\right)\dif s\right)^{\frac{2}{\beta}}\no\\
&\leq C/R^{\beta-2}\to 0\mbox{ uniformly in $n,\eps$ as $R\to\infty$}.\label{LQ1}
\end{align}
For $J^{n,2}_{R,\eps}(t)$, we can not directly use the discretized Krylov estimate to conclude
\begin{align}\label{Lim0}
\lim_{\eps\to 0}\sup_n J^{n,2}_{R,\eps}(t)=0,\ \ \forall t,R>0,
\end{align}
because the Euler scheme \eqref{EU21} has unbounded coefficients.
We need to cutoff the coefficients. Let $\chi_R(x)$ be a nonnegative smooth cutoff function with $\chi_R(x)=1$ for $|x|<R$ and $\chi_R(x)=0$ for $|x|>R+1$.
Define 
$$
b^{n,R}_t(x):=b_t(x,\mu_{ X^n_{t_n}})\chi_R(x), \ \sigma^{n,R}_t(x):=\sigma_t(x\chi_R(x),\mu_{ X^n_{s_n}}).
$$
Let $(X^{n,R},\bar X^{n,R})$ solve the following equation in $\mR^{2d}$ (no coupling):
\begin{align*}
\left\{
\begin{aligned}
\dif  X^{n,R}_t=b^{n,R}_{t_n}\big(X^{n,R}_{t_n}\big)\dif t+\sigma^{n,R}_{t_n}\big(X^{n,R}_{t_n}\big)\dif W_t,\ 
X^{n,R}_0=X^n_0,\\
\dif  \bar X^{n,R}_t=b^{n,R}_{t_n}\big(\bar X^{n,R}_{t_n}\big)\dif t+\sigma^{n,R}_{t_n}\big(\bar X^{n,R}_{t_n}\big)\dif\bar W_t,\ \bar X^{n,R}_0=\bar X^n_0,
\end{aligned}
\right.
\end{align*}
where $(W,X^n_0)$ and $(\bar W, \bar X^n_0)$ are the same as in \eqref{EU21}.
From the construction, one sees that
\begin{align}\label{EU3}
(X^{n}_t,\bar X^{n}_t)=(X^{n,R}_t,\bar X^{n,R}_t),\ \ t<\tau^n_R.
\end{align}
Moreover, it is easy to see that $(X^{n,R}, \bar X^{n,R})$ is a discretized $\mR^{2d}$-valued It\^o process with coefficients satisfying the assumptions in Theorem \ref{Th0} uniformly in $n$.
Thus, for fixed $R>0$ and any $p>2d+1$, by \eqref{EU3} and \eqref{EY1} we have
\begin{align*}
J^{n,2}_{R,\eps}(t)&=\int^t_0\mE\Big(\1_{t<\tau^n_R}\|\bar\sigma_{s_n}(  X^{n,R} _{s_n},\bar X^{n,R}_{s_n})
-\bar\sigma^\eps_{s_n}( X^{n,R}_{s_n},\bar X^{n,R}_{s_n})\|^2\Big)\dif s\\
&\leq\int^t_0\mE\Big(\1_{|X^{n,R} _{s_n}|\vee|\bar X^{n,R}_{s_n}|<R}
\|\bar\sigma_{s_n}(  X^{n,R} _{s_n},\bar X^{n,R}_{s_n})-\bar\sigma^\eps_{s_n}( X^{n,R}_{s_n},\bar X^{n,R}_{s_n})\|^2\Big)\dif s\\
&=\int^{h_n}_0\mE\Big(\1_{|X^{n} _0|\vee|\bar X^{n}_0|<R}
\|\bar\sigma_0(  X^{n}_0,\bar X^{n}_0)-\bar\sigma^\eps_0( X^{n}_0,\bar X^{n}_0)\|^2\Big)\dif s\\
&+\int^t_{h_n}\mE\Big(\1_{|X^{n,R} _{s_n}|\vee|\bar X^{n,R}_{s_n}|<R}
\|\bar\sigma_{s_n}(  X^{n,R} _{s_n},\bar X^{n,R}_{s_n})-\bar\sigma^\eps_{s_n}( X^{n,R}_{s_n},\bar X^{n,R}_{s_n})\|^2\Big)\dif s\\
&\lesssim h_n+\left(\int^t_0\|\bar\sigma_{s_n}( \cdot,\cdot)-\bar\sigma^\eps_{s_n}(\cdot,\cdot)\|^{2p}_{L^{2p}(B_R\times B_R)}\dif s\right)^{1/p}\\
&\lesssim h_n+ \left(\int^t_0\|\bar\sigma_s( \cdot,\cdot)-\bar\sigma^\eps_s(\cdot,\cdot)\|_{L^{2p}(B_R\times B_R)}^{2p}\dif s\right)^{1/p}\\
&\quad+\left(\int^t_0\|\bar\sigma_s(\cdot,\cdot)-\bar\sigma_{s_n}(\cdot,\cdot)\|_{L^{2p}(B_{R+1}\times B_{R+1})}^{2p}\dif s\right)^{1/p}\\
&=:h_n+I^{R}_\eps(t)+K^{R}_n(t),
\end{align*}
where $h_n\downarrow 0$ as $n\to\infty$, and the constants contained in the above $\lesssim$ may depend on $R$.
By the dominated convergence theorem and the continuity of $t\mapsto \sigma_t(x,y)$, we have
$$
 \lim_{\eps\to\infty}I^{R}_\eps(t)=0,\ \ \lim_{n\to\infty}K^{R}_n(t)=0,
$$
which in turn implies the limit \eqref{Lim0}, and so \eqref{KQ1}. 
Thus we complete the proof.
\end{proof}

\begin{proof}[Proof of (i) of Theorem \ref{Th12}]
Using the above lemma and taking limits for both sides of \eqref{KK1}, one finds that $(\tilde X,\tilde W)$ solves the following SDE:
\begin{align}\label{KK2}
\tilde X _t=\tilde X_0+\int^t_0b_s\left(\tilde X_{s},\mu_{\tilde X_{s}}\right)\dif s
+\int^t_0\sigma_s\left(\tilde X_{s},\mu_{\tilde X_{s}}\right)\dif \tilde W_s.
\end{align}
Since the weak uniqueness holds for DDSDE \eqref{SDE0}, any weak solutions have the same distribution. Hence, the whole Euler approximation $X^h$
weakly converges to the unique weak solution $X$ in distribution. As for \eqref{WE}, it follows by Krylov's estimate \eqref{EY1}.
\end{proof}

In order to show (ii) of Theorem \ref{Th12}, we need the following important observation due to \cite[Lemma 1.1]{Gy-Kr}, which has the root of Yamada-Watanabe's theorem.
\bl\label{Le34}
Let $(Z_h)_{h\in(0,1)}$ be a family of random elements in a Polish space $(E,\rho)$. Then $Z_h$ converges in probability to an $E$-valued random element as $h\to 0$ 
if and only if for every pair of subsequences $(Z_{h_n}, Z_{\ell_n})_{n\in\mN}$, there exists a subsubsequence  $(Z_{h_{n(k)}}, Z_{\ell_{n(k)}})_{k\in\mN}$ 
converging in distribution to a random element in $E\times E$, which supports on the diagonal $\{(x,y)\in E\times E: x=y\}$.
\el
\begin{proof}
We use a contradiction method. Suppose that $Z_h$ does not converge in probability. Then there is an $\eps>0$ such that for any $\delta>0$, there are
$h_\delta$ and $\ell_\delta$ less than $\delta$ such that
$$
\mP\big\{\rho(Z_{h_\delta}, Z_{\ell_\delta})>\eps\big\}\geq\eps.
$$
Thus we can choose two subsequence $Z_{h_n}$ and $Z_{\ell_n}$ such that 
\begin{align}\label{KA1}
\inf_{n\in\mN}\mP\big\{\rho(Z_{h_n}, Z_{\ell_n})>\eps\big\}\geq\eps.
\end{align}
By the assumption, there is a subsubsequence $(Z_{h_{n(k)}}, Z_{\ell_{n(k)}})_{k\in\mN}$ such that
$$
\lim_{k\to\infty}\mE\left(\rho(Z_{h_{n(k)}}, Z_{\ell_{n(k)}})\wedge 1\right)=0.
$$
Clearly, this is contradict with \eqref{KA1}. By the completeness of $(E,\rho)$, we complete the proof.
\end{proof}

Now we are in a position to give

\begin{proof}[Proof of (ii) of Theorem \ref{Th12}]
Let $X^{h_n}$ and $X^{\ell_n}$ be two subsequences of $X^h$. Clearly, by Lemma \ref{Le31},
the law  of $(X^{h_n},X^{\ell_n}, W)_{ n\in\mN}$  in $\mC\times\mC\times\mC$ is tight.
As above, by Skorokhod's embedding theorem, there exist subsequences $n(k)$, a probability space $(\tilde\Omega,\tilde\sF,\tilde\mP)$,
carrying stochastic processes $(\tilde X^{h_{n(k)}},\hat X^{\ell_{n(k)}},\tilde  W^k)$ and $(\tilde X,\hat X, \tilde W)$ such that
$$
\left(\tilde X^{h_{n(k)}},\hat X^{\ell_{n(k)}},\tilde  W^k\right)\stackrel{k\to\infty}{\to}\left (\tilde X,\hat X, \tilde W\right)\ \ \tilde\mP-a.s.
$$
and for each $k\in\mN$,
$$
\tilde\mP\circ\left(\tilde X^{h_{n(k)}},\hat X^{\ell_{n(k)}},\tilde  W^k\right)^{-1}=\mP\circ\left(X^{h_{n(k)}},X^{\ell_{n(k)}},W\right)^{-1}.
$$
As in showing \eqref{KK2}, one sees that $(\tilde X,\tilde W)$ and $(\hat X,\tilde W)$ are two solutions of DDSDE \eqref{SDE0} 
defined on the same probability space with the same initial values $\tilde X_0=\hat X_0$. The latter point is due to
$$
\tilde\mP(\tilde X_0=\hat X_0)\geq \lim_{k\to\infty}\tilde\mP\left(\tilde X^{h_{n(k)}}_0=\hat X^{\ell_{n(k)}}_0\right)
=\lim_{k\to\infty}\mP\left(X^{h_{n(k)}}_0=X^{\ell_{n(k)}}_0\right)=1.
$$
By the pathwise uniqueness, we obtain $\tilde X=\hat X$. Thus by Lemma \ref{Le34}, we conclude that $X^h$ converges in probability to a random elelment $X$ in $\mC$
as $h\downarrow 0$.
Using Lemma \ref{Le32}, one sees that $X$ is a solution of  DDSDE \eqref{SDE0}. Moreover, the convergence \eqref{WE0} follows by \eqref{KJ1}
and the dominated convergence theorem.
\end{proof}

\section{Propagation of chaos: Proof of Theorem \ref{Th13}}

In this section we use induction  to prove Theorem \ref{Th13}.
First of all, we prepare several lemmas. The following lemma is the same as in Lemma \ref{Le31}. We omit the details.
\bl\label{Le41}
Let $\beta\geq 2$. Under {\bf (H$'_\beta$)}, for any $T>0$, there is a constant $C>0$ such that for all $N\in\mN$ and $j=1,\cdots, N$,
\begin{align}\label{KJ11}
\mE\left(\sup_{t\in[0,T]}|X_t^{N,j}|^\beta\right)\leq C(1+\mE|\xi_j|^\beta).
\end{align}
\el
\bl
Let $\beta\geq 2$. Under {\bf (H$'_\beta$)}, for any $T>0$, there is a constant $C>0$ such that for all $N\in\mN$ and $j=1,\cdots,N$,
\begin{align}\label{UY1}
\mE\left|\frac{1}{N}\sum_{i=1}^N\Big(b_t(\bar X^{j}_t, \mu_{\bar X^j_t})-\bar b_t(\bar X^{j}_t, \bar X^i_t)\Big)\right|^2\leq C/N.
\end{align}
\el
\begin{proof}
Notice that 
$$
\mu_{\bar X^1_t}=\cdots=\mu_{\bar X^N_t}.
$$
For simplicity, if we define
$$
\hat b_t(x,y):=b_t(x, \mu_{\bar X^1_t})-\bar b_t(x, y),
$$
then the left hand side of \eqref{UY1} denoted by $\sI$ can be written as
\begin{align*}
\sI=\mE\left|\frac{1}{N}\sum_{i=1}^N\hat b_t(\bar X^{j}_t,\bar X^i_t)\right|^2
=\frac{1}{N^2}\sum_{i,k=1}^N\mE\<\hat b_t(\bar X^{j}_t,\bar X^i_t),\hat b_t(\bar X^{j}_t,\bar X^k_t)\>.
\end{align*}
For $i\not=j\not=k$, since $\bar X^{i}_t, \bar X^{j}_t, \bar X^{k}_t$ are independent, we have
$$
\mE\<\hat b_t(\bar X^{j}_t,\bar X^i_t),\hat b_t(\bar X^{j}_t,\bar X^k_t)\>=0.
$$
Therefore, by Lemma \ref{Le41} and the linear growth of $\hat b$, we get
$$
\sI\leq \frac{2}{N^2}\sum_{i=1}^N\mE|\hat b_t(\bar X^j_t,\bar X^i_t)|^2
\leq \frac{C}{N^2}\sum_{i=1}^N\Big(1+\mE|\bar X^j_t|^2+\mE|\bar X^i_t|^2\Big)\leq  C/N.
$$
The proof is complete.
\end{proof}
\bl\label{Le43}
Let $m\in\mN$ and $f:\mR_+\times\mR^d\times\mR^d\to \mR^m$ be a locally bounded measurable function with
$$
|f_t(x,y)|\leq c(1+|x|+|y|).
$$
Let $f^\eps_t(x,y):=f_t(\cdot)*\varrho_\eps(x,y)$ be the mollifying approximation. Define
\begin{align}
\sA^f_{N,\eps}(t)&:=\sup_{j=1,\cdots, N}\mE\int^t_0|f_{s_h}-f^\eps_{s_h}|^2(X^{N,j}_{s_h}, \mu^{N}_{s_h})\dif s,\label{DF1}\\
\bar\sA^f_{N,\eps}(t)&:=\sup_{j=1,\cdots, N}\mE\int^t_0|f_{s_h}-f^\eps_{s_h}|^2(\bar X^{j}_{s_h}, \mu_{\bar X^j_{s_h}})\dif s.\label{DF2}
\end{align}
Let $\beta>2$ and $\nu\in\cP_\beta(\mR^d)$. Suppose that {\bf (H$'_\beta$)} holds 
and  the initial law $\nu$ has a density $\phi\in L^q_{loc}(\mR^d)$ for some $q>1$. Then we have
$$
\lim_{\eps\to 0}\varlimsup_{N\to\infty} \sA^f_{N,\eps}(t)=0,\ \ \lim_{\eps\to 0}\varlimsup_{N\to\infty}\bar\sA^f_{N,\eps}(t)=0.
$$
\el
\begin{proof}
We only prove the first limit. For simplicity, we write
$$
F^\eps_t(x,y):=|f_t(x,y)-f^\eps_t(x,y)|^2.
$$
Without loss of generality we assume $t>h$. Notice that
\begin{align*}
\mE\int^t_0F^\eps_{s_h}(X^{N,j}_{s_h}, \mu^{N}_{s_h})\dif s&=
h\mE F^\eps_0(X^{N,j}_0, \mu^{N}_0)+
\mE\int^t_hF^\eps_{s_h}(X^{N,j}_{s_h}, \mu^{N}_{s_h})\dif s=:I^{(1)}_{N,\eps}+I^{(2)}_{N,\eps}.
\end{align*}
For $I^{(1)}_{N,\eps}$, by the assumption we have
\begin{align*}
I^{(1)}_{N,\eps}&=\frac{h}{N}\sum_{i=1}^N\mE F^\eps_0(X^{N,j}_0,X^{N,i}_0)=\frac{h}{N}\int_{\mR^d}F^\eps_0(x,x)\phi(x)\dif x\\
&+\frac{h(N-1)}{N}\int_{\mR^{2d}}F^\eps_0(x,y)\phi(x)\phi(y)\dif x\dif y=:I^{(11)}_{N,\eps}+I^{(12)}_{N,\eps}.
\end{align*}
For $I^{(11)}_{N,\eps}$, we clearly have
$$
\lim_{N\to\infty}\sup_\eps I^{(11)}_{N,\eps}\leq \lim_{N\to\infty}\frac{C}{N}\int_{\mR^d}(1+|x|^2)\phi(x)\dif x=0.
$$
For $I^{(12)}_{N,\eps}$, if we define $B_R:=\{(x,y)\in\mR^d\times\mR^d: |x|<R,|y|<R\}$ for $R>0$, then by H\"older's inequality and $\phi\in L^q_{loc}(\mR^d)$,
\begin{align*}
I^{(12)}_{N,\eps}&\leq\int_{B_R}F^\eps_0(x,y)\phi(x)\phi(y)\dif x\dif y+\int_{B_R^c}F^\eps_0(x,y)\phi(x)\phi(y)\dif x\dif y\\
&\lesssim\left(\int_{B_R}|F^\eps_0(x,y)|^{\frac{q}{q-1}}\dif x\dif y\right)^{\frac{q-1}{q}}\left(\int_{B_R}|\phi(x)\phi(y)|^q\dif x\dif y\right)^{\frac{1}{q}}\\
&\quad+\int_{B_R^c}(1+|x|^2+|y|^2)\phi(x)\phi(y)\dif x\dif y\\
&\lesssim\left(\int_{B_R}|F^\eps_0(x,y)|^{\frac{q}{q-1}}\dif x\dif y\right)^{\frac{q-1}{q}}\left(\int_{|x|<R}|\phi(x)|^q\dif x\right)^{\frac{2}{q}}\\
&\quad+\frac{1}{R^{\beta-2}}\int_{\mR^{2d}}(1+|x|^\beta+|y|^\beta)\phi(x)\phi(y)\dif x\dif y,
\end{align*}
where the constant $C$ contained in $\lesssim$ is independent of $N,R, \eps$.
By the dominated convergence theorem and first letting $\eps\to 0$ and then $R\to\infty$, we get
$$
\lim_{\eps\to 0}\sup_N I^{(12)}_{N,\eps}=0.
$$
Next we treat $I^{(2)}_{N,\eps}$ and write
\begin{align*}
I^{(2)}_{N,\eps}&=\frac{1}{N}\mE\int^t_hF^\eps_{s_h}(X^{N,1}_{s_h},X^{N,1}_{s_h})\dif s
+\frac{1}{N}\sum_{i=2}^N\mE\int^t_hF^\eps_{s_h}(X^{N,1}_{s_h},X^{N,i}_{s_h})\dif s=:I^{(21)}_{N,\eps}+I^{(22)}_{N,\eps}.
\end{align*}
For $I^{(21)}_{N,\eps}$, by \eqref{KJ11} we have
$$
\lim_{N\to\infty}\sup_\eps I^{(21)}_{N,\eps}\leq\lim_{N\to\infty}\frac{C}{N}\mE\int^t_h(1+|X^{N,1}_{s_h}|^2)\dif s=0.
$$
For $I^{(22)}_{N,\eps}$, using the discretized Krylov estimate and the same argument as in showing \eqref{Lim1}, we also have
$$
\lim_{\eps\to 0}\sup_N I^{(22)}_{N,\eps}=0.
$$
Combining the above limits, we complete the proof.
\end{proof}

Now we can give 
\begin{proof}[Proof of Theorem \ref{Th13}]
By equations \eqref{SDE99}, \eqref{SDE98} and Burkholder's inequality we have
\begin{align*}
\mE\left(\sup_{s\in[0,t]}|X^{N,j}_s-\bar X^j_s|^2\right)
&\lesssim\mE\int^t_0\Big|b_{s_h}\big(X^{N,j}_{s_h}, \mu^{N}_{s_h}\big)-b_{s_h}\big(\bar X^{j}_{s_h}, \mu_{\bar X^j_{s_h}}\big)\Big|^2\dif s\\
&+\mE\int^t_0\Big|\sigma_{s_h}\big(X^{N,j}_{s_h}, \mu^{N}_{s_h}\big)-\sigma_{s_h}\big(\bar X^{j}_{s_h}, \mu_{\bar X^j_{s_h}}\big)\Big|^2\dif s\\
&\lesssim\mE\int^t_0\Big|b^\eps_{s_h}\big(X^{N,j}_{s_h}, \mu^{N}_{s_h}\big)-b^\eps_{s_h}\big(\bar X^{j}_{s_h}, \mu_{\bar X^j_{s_h}}\big)\Big|^2\dif s\\
&+\mE\int^t_0\Big|\sigma^\eps_{s_h}\big(X^{N,j}_{s_h}, \mu^{N}_{s_h}\big)-\sigma^\eps_{s_h}\big(\bar X^{j}_{s_h}, \mu_{\bar X^j_{s_h}}\big)\Big|^2\dif s\\
&+\sA^b_{N,\eps}(t)+\bar\sA^b_{N,\eps}(t)+\sA^\sigma_{N,\eps}(t)+\bar\sA^\sigma_{N,\eps}(t),
\end{align*}
where for $f=b$ and $\sigma$, $f^\eps_t(x,y):=f_t*\varrho_\eps(x,y)$ is the mollifying approximation, 
and $\sA^f_{N,\eps}(t)$ and $\bar\sA^f_{N,\eps}(t)$ are defined by \eqref{DF1} and \eqref{DF2}, respectively.
Now let $\chi_R(x,y)$ be a smooth function with 
$$
\chi_R(x,y)=1,\ |x|\vee|y|<R,\ \chi_R(x,y)=0,\ \ |x|\vee|y|>R+1.
$$
For $f=b$ or $\sigma$, by definition, we make the following decomposition:
\begin{align*}
&\mE\int^t_0\Big|f^\eps_{s_h}\big(X^{N,j}_{s_h}, \mu^{N}_{s_h}\big)-f^\eps_{s_h}\big(\bar X^{j}_{s_h}, \mu_{\bar X^j_{s_h}}\big)\Big|^2\dif s\\
&=\mE\int^t_0\Big|\frac{1}{N}\sum_{i=1}^N\left(f^\eps_{s_h}\big(X^{N,j}_{s_h}, X^{N,i}_{s_h}\big)
-f^\eps_{s_h}\big(\bar X^{j}_{s_h}, \bar X^i_{s_h}\big)\right)\Big|^2\dif t\\
&\leq\frac{1}{N}\sum_{i=1}^N\mE\int^t_0\left|f^{\eps}_{s_h}\big(X^{N,j}_{s_h}, X^{N,i}_{s_h}\big)
-f^{\eps}_{s_h}\big(\bar X^{j}_{s_h}, \bar X^i_{s_h}\big)\right|^2\dif t\\
&\leq\frac{2}{N}\sum_{i=1}^N\mE\int^t_0\left|f^{\eps,R}_{s_h}\big(X^{N,j}_{s_h}, X^{N,i}_{s_h}\big)
-f^{\eps,R}_{s_h}\big(\bar X^{j}_{s_h}, \bar X^i_{s_h}\big)\right|^2\dif s\\
&\quad+\frac{4}{N}\sum_{i=1}^N\mE\int^t_0|f^\eps_{s_h}-f^{\eps,R}_{s_h}|^2\big(X^{N,j}_{s_h}, X^{N,i}_{s_h}\big)\dif s\\
&\quad+\frac{4}{N}\sum_{i=1}^N\mE\int^t_0|f^\eps_{s_h}-f^{\eps,R}_{s_h}|^2\big(\bar X^{j}_{s_h}, \bar X^i_{s_h}\big)\dif s\\
&=:\frac{1}{N}\sum_{i=1}^N\left(\sI^{i,j}_{R,1}+\sI^{i,j}_{R,2}+\sI^{i,j}_{R,3}\right),
\end{align*}
where
$$
f^{\eps,R}_{s_h}(x,y)=f^{\eps}_{s_h}(x,y)\chi_R(x,y).
$$
Since $f$ is linear growth, we have
$$
|f^\eps_{s_h}-f^{\eps,R}_{s_h}|(x,y)\leq C (\1_{|x|>R}+\1_{|y|>R})(1+|x|+|y|).
$$
Thus, by H\"older's inequality and \eqref{KJ11}, we have
\begin{align*}
\sI^{i,j}_{R,2}&\lesssim \mE\int^t_0\left(\1_{|X^{N,j}_{s_h}|>R}+\1_{|X^{N,i}_{s_h}|>R}\right)\big(1+|X^{N,j}_{s_h}|^2+|X^{N,i}_{s_h}|^2\big)\dif s\\
&\lesssim \int^t_0\left(\mP(|X^{N,j}_{s_h}|>R)+\mP(|X^{N,i}_{s_h}|>R)\right)^{\frac{\beta-2}{\beta}}\dif s\\
&\lesssim \frac{1}{R^{\beta-2}}\int^t_0\left(\mE|X^{N,j}_{s_h}|^\beta+\mE|X^{N,i}_{s_h}|^\beta\right)^{\frac{\beta-2}{\beta}}\dif s
\leq \frac{C}{R^{\beta-2}},
\end{align*}
where $C$ is independent of $N, i,j$. Similarly, we also have
$$
\sup_{i,j}\sI^{i,j}_{R,3}\leq \frac{C}{R^{\beta-2}}.
$$
For $\sI^{i,j}_{R,1}$, since $(x,y)\mapsto f^{\eps,R}_{s_h}(x,y)$ is Lipschitz continuous, we have
$$
\sI^{i,j}_{R,1}\leq C_{\eps,R}\mE\int^t_0\Big(|X^{N,j}_{s_h}-\bar X^{j}_{s_h}|^2+|X^{N,i}_{s_h}-\bar X^{i}_{s_h}|^2\Big)\dif s.
$$
Combining the above calculations we obtain that for all $t\in[0,T]$,
\begin{align}
&\sup_{j=1,\cdots,N}\mE\left(\sup_{s\in[0,t]}|X^{N,j}_s-\bar X^j_s|^2\right)
\leq C_{\eps,R}\sup_{j=1,\cdots,N}\int^t_0\mE |X^{N,j}_{s_h}-\bar X^{j}_{s_h}|^2\dif s\no\\
&\qquad\qquad +C/R^{\beta-2}+\sA^b_{N,\eps}(t)+\bar\sA^b_{N,\eps}(t)+\sA^\sigma_{N,\eps}(t)+\bar\sA^\sigma_{N,\eps}(t).\label{DF4}
\end{align}
Here we can not use Gronwall's inequality to derive the result. We shall use the induction method to show \eqref{DF3}.
First of all, we clearly have
\begin{align*}
\mE |X^{N,j}_0-\bar X^{j}_0|^2=0.
\end{align*}
Suppose that we have shown that for some $k\in\mN_0$,
$$
\lim_{N\to\infty}\sup_{j=1,\cdots,N}\mE\left(\sup_{s\in[0,kh]}|X^{N,j}_s-\bar X^j_s|^2\right)=0.
$$
Then for $t=(k+1)h$, by \eqref{DF4} we have
\begin{align*}
&\sup_{j=1,\cdots,N}\mE\left(\sup_{s\in[0,t]}|X^{N,j}_s-\bar X^j_s|^2\right)
\leq C_{\eps,R}h \sup_{j=1,\cdots,N}\sum_{m=0}^k \mE |X^{N,j}_{mh}-\bar X^{j}_{mh}|^2\no\\
&\qquad\qquad +C/R^{\beta-2}+\sA^b_{N,\eps}(t)+\bar\sA^b_{N,\eps}(t)+\sA^\sigma_{N,\eps}(t)+\bar\sA^\sigma_{N,\eps}(t).\label{DF4}
\end{align*}
Firstly letting $N\to\infty$ and then $R\to\infty$ and $\eps\to 0$, by Lemma \ref{Le43} and the induction hypothesis, we obtain \eqref{DF3}
for $t=(k+1)h$. The proof is complete.
\end{proof}

\medskip

\begin{center}
\bf Acknowledgement
\end{center}

{The author thanks Jianhai Bao and Zimo Hao for their useful discussions and suggestions.}

\end{document}